\documentclass[11pt]{amsart}
\usepackage{graphicx}

\usepackage{amsmath}
\usepackage{amsthm}
\usepackage{amssymb}
\usepackage{enumerate}
\usepackage{cleveref}

\newtheorem{lem}{Lemma}[section]

\newtheorem{prop}[lem]{Proposition}
\newtheorem{coro}[lem]{Corollary}
\newtheorem{thm}[lem]{Theorem}
\newtheorem{conj}[lem]{Conjecture}
\newtheorem{de}[lem]{Definition}

\newtheorem{pozn}[lem]{Remark}
\newtheorem{ex}[lem]{Example}
\newtheorem{klejm}[lem]{Claim}
\newtheorem{obs}[lem]{Observation}
\def\pf{\begin{proof}}
\def\pfk{\end{proof}}

\begin{document}
\title{On Brlek-Reutenauer conjecture}
\maketitle

\begin{center}
\author{L{\!'}. Balkov\'a\footnote{e-mail: lubomira.balkova@fjfi.cvut.cz}, E. Pelantov\'a, \v S. Starosta\\[2mm]
{\normalsize Department of Mathematics FNSPE, Czech Technical
University in Prague\\ and Doppler Institute for Mathematical
Physics and
Applied Mathematics}\\
{\normalsize Trojanova 13, 120 00 Praha 2, Czech Republic}}

\end{center}

\begin{abstract} Brlek and Reutenauer conjectured that any infinite word ${\mathbf u}$ with language closed under
reversal satisfies the equality $2D({\mathbf u}) =
\sum_{n=0}^{+\infty}T_{\mathbf u}(n)$ in which  $D({\mathbf u})$
denotes the defect of $\mathbf u$ and $T_{\mathbf u}(n)$ denotes
$\mathcal{C}_{\mathbf u}(n+1)-\mathcal{C}_{\mathbf u}(n)  +2 -
\mathcal{P}_{\mathbf u}(n+1) - \mathcal{P}_{\mathbf u}(n)$, where
$\mathcal{C}_{\mathbf u}$  and $ \mathcal{P}_{\mathbf u}$ are
the factor and palindromic complexity of ${\mathbf u}$,
respectively. Brlek and Reutenauer verified their conjecture for
periodic infinite words. We prove the conjecture for uniformly
recurrent words. Moreover, we summarize results and some open
problems related to defect, which may be useful for the proof of
Brlek-Reutenauer Conjecture in full generality.
\end{abstract}

\section{Introduction}

There have been recently quite a lot of papers devoted to palindromes in infinite words.
Droubay, Justin, and Pirillo determined in~\cite{DrJuPi} the upper bound on the number of distinct
palindromes occurring in a~finite word -- a finite word $w$ contains at most $|w|+1$ different palindromes, where
$|w|$ denotes the length of $w$. The difference between the utmost number $|w|+1$ and the actual number of palindromes in $w$
is called the defect of $w$ and is usually denoted by $D(w)$. An infinite word $\mathbf u$ whose factors
have all zero defect was baptized rich or full.
In~\cite{BaMaPe}, Bal\'a\v zi, Mas\'akov\'a, and Pelantov\'a proved for infinite words with language closed under reversal
an inequality relating the palindromic and factor complexity of an infinite word $\mathbf u$ denoted ${\mathcal P}_{\mathbf u}$ and ${\mathcal C}_{\mathbf u}$, respectively.
For such infinite words, it holds
\begin{equation}\label{BalaziNerovnostUvod}
{\mathcal C}_{\mathbf u}(n+1)-{\mathcal C}_{\mathbf u}(n)+2-{\mathcal P}_{\mathbf u}(n)-{\mathcal P}_{\mathbf u}(n+1)\geq 0 \quad \text{for all $n \in \mathbb N$}.
\end{equation}
In~\cite{BuLuGlZa}, Bucci, De Luca, Glen, and Zamboni showed that rich words with language closed under reversal can be characterized by the equality in~\eqref{BalaziNerovnostUvod}.
Brlek, Hamel, Nivat, and Reutenauer in~\cite{BHNR} defined
the defect $D({\mathbf u})$ of an infinite word ${\mathbf u}$ as the maximum
of defects of all its factors and they studied its value for periodic words.

Recently, in~\cite{BaPeSt}, the authors of this paper have proven
that for a~uniformly recurrent word $\mathbf u$, its defect $D(\mathbf u)$ is finite if and only if the equality in~\eqref{BalaziNerovnostUvod} is attained for all but a~finite number of indices~$n$.

Despite the fact that numerous researchers study palindromes, only
recently Brlek and Reute\-nauer have noticed that the value of
defect is closely tied with the expression on the left-hand side
of~\eqref{BalaziNerovnostUvod} - let us denote it by $T_{\mathbf
u}(n)$. They have shown that for periodic infinite words with
language closed under reversal, it holds $2D(\mathbf
u)=\sum_{n=0}^{+\infty}T_{\mathbf u}(n)$. Their conjecture says
that the same equation holds for all infinite words with language
closed under reversal.

In this paper we will prove that Brlek-Reutenauer Conjecture is true for uniformly recurrent words and in the last chapter
we will discuss some aspects concerning the conjecture for infinite words that are not uniformly recurrent.

\section{Preliminaries}
By $\mathcal{A}$ we denote a~finite set of symbols called
{\em letters}; the set $\mathcal{A}$ is therefore called an {\em
alphabet}. A finite string $w=w_0w_1\ldots w_{n-1}$ of letters from
$\mathcal{A}$ is said to be a~{\em finite word}, its length is
denoted by $|w| = n$. Finite words over  $\mathcal{A}$ together
with the operation of concatenation and the empty word $\epsilon$
as the neutral element form a~free monoid $\mathcal{A}^*$. The map
$$w=w_0w_1\ldots w_{n-1} \quad \mapsto \quad \overline{w} =
w_{n-1}w_{n-2}\ldots w_{0}$$ is a bijection on $\mathcal{A}^*$,
the word $\overline{w}$ is called the {\em reversal} or the {\em
mirror image} of $w$. A~word $w$ which coincides with its mirror
image is a~{\em palindrome}.

Under an {\em infinite word} we understand an infinite string
${\mathbf u}=u_0u_1u_2\ldots $ of letters from $\mathcal{A}$.
A~finite word $w$ is a~{\em factor} of a~word $v$ (finite or
infinite) if there exist words $p$ and $s$ such that $v= pws$. If
$p = \epsilon$, then $w$ is said to be a~{\em prefix} of $v$, if
$s = \epsilon$, then $w$ is a~{\em suffix} of~$v$.

The {\em language} $\mathcal{L}({\mathbf u})$ of an infinite word
${\mathbf u}$ is the set of all its factors. Factors of $\mathbf
u$ of length $n$ form the set denoted by $\mathcal{L}_n({\mathbf
u})$. We say that the language
$\mathcal{L}({\mathbf u})$ is {\em closed under reversal} if
$\mathcal{L}({\mathbf u})$ contains with every factor $w$ also its
reversal $\overline{w}$.

For any factor $w\in \mathcal{L}({\mathbf u})$, there exists an
index $i$ such that $w$ is a prefix  of  the infinite word
$u_iu_{i+1}u_{i+2} \ldots$. Such an index is called an {\em
occurrence} of $w$ in ${\mathbf u}$. If each factor of $\mathbf u$
has infinitely many occurrences in ${\mathbf u}$, the infinite
word $\mathbf u$ is said to be {\em recurrent}. It is easy to see
that if the language of ${\mathbf u}$ is closed under reversal,
then ${\mathbf u}$ is recurrent (a~proof can be found
in~\cite{GlJuWiZa}). For a~recurrent infinite word ${\mathbf u}$,
we may define the notion of a~{\em complete return word} of any $w
\in\mathcal{L}({\mathbf u})$. It is a~factor $v\in
\mathcal{L}({\mathbf u})$ such that $w$ is a prefix and a suffix
of $v$ and $w$ occurs in $v$ exactly twice. Under a~{\em return
word} of a factor $w$ is usually understood a word $q \in
\mathcal{L}({\mathbf u})$  such that $qw$ is a complete return
word of $w$. If any factor $w \in \mathcal{L}({\mathbf u})$ has
only finitely many return words, then the infinite word ${\mathbf
u}$ is called {\em uniformly recurrent}. If ${\mathbf u}$ is
a~uniformly recurrent word, we can find for any $n \in \mathbb{N}$
a~number $R$ such that any factor of ${\mathbf u}$ which is longer
than $R$ contains already all factors of ${\mathbf u}$
of length $n$.

The {\em factor complexity} of an infinite word ${\mathbf u}$ is
a~map $\mathcal{C}_{\bf u}: \mathbb{N} \mapsto \mathbb{N}$ defined
by the prescription $\mathcal{C}_{\bf u}(n):=\#
\mathcal{L}_n({\mathbf u})$. To determine the first difference of
the factor complexity, one has to count the possible {extensions}
of factors of length $n$. A~{\em right extension} of $w \in
\mathcal{L}({\mathbf u})$ is any letter $a\in \mathcal{A}$ such
that $w a\in \mathcal{L}({\mathbf u})$. Of course, any factor of
${\mathbf u}$ has at least one right extension. A~factor $w$ is
called {\em right special} if $w$ has at least two right
extensions. Similarly, one can define a~{\em left extension} and
a~{\em left special} factor. We will deal mainly with recurrent
infinite words ${\mathbf u}$. In such a~case, any factor of $\mathbf
u$ has at least one left extension.

\medskip

The {{\em defect} $D(w)$ of a~finite word $w$ is the difference
between the utmost number of palindromes $|w|+1$ and the actual
number of palindromes contained in $w$. Finite words with zero
defects -- called {\em rich} or {\em full} words -- can be viewed
as the most saturated by palindromes. This definition may be
extended to infinite words as follows.}

\begin{de} An infinite word ${\mathbf u} = u_0u_1u_2 \ldots$ is called
{\it rich} or {\it full}, if for any index $n\in \mathbb{N}$  the
prefix $u_0u_1u_2 \ldots u_{n-1}$ of length $n$ contains exactly
$n+1$ different palindromes.
\end{de}

Let us remark that not only all prefixes of rich words are rich,
but also all factors are rich. A~result from~\cite{DrJuPi} will
provide us with a~handful tool which helps to evaluate the defect
of a factor.

\begin{prop}[\cite{DrJuPi}]\label{lps} A finite or infinite word $\mathbf u$ is rich if and only if
the longest palindromic suffix of $w$ occurs exactly once in $w$
for any prefix $w$ of $\mathbf u$.
\end{prop}

In accordance  with the terminology introduced in \cite{DrJuPi},
the factor with a~unique occurrence in another factor is called
{\em unioccurrent}. From the proof of the previous proposition
directly follows the next corollary.

\begin{coro}\label{prodlouzeni} The defect $D(w)$ of a~finite word $w$ is equal to the
number of prefixes $w'$ of $w$, for which the longest palindromic
suffix of $w'$ is not unioccurrent in $w'$. In other words, if $b$
is a~letter and $w$ a~finite word, then $D(wb)=D(w)+\delta$, where
$\delta = 0$ if the longest palindromic suffix of $wb$ occurs
exactly once in $wb$ and $\delta = 1$ otherwise.
\end{coro}
This corollary implies that $D(v) \geq D(w)$ whenever $w$ is
a~factor of $v$. It enables to give a~reasonable definition of the
defect of an infinite word (see~\cite{BHNR}).

\begin{de}  The defect of  an  infinite word  $\mathbf{u}$ is the  number (finite or infinite)
$$D({\mathbf{u}}) = \sup \{ D(w) \bigm | w \ \text{is a~prefix of $\mathbf u$}\}\,.$$
\end{de}
Let us point out several facts concerning defects that are easy to
prove:
\begin{enumerate}
\item If we consider all factors of a~finite or an infinite word
$\mathbf u$, we obtain the same defect, i.e.,
$$D({\mathbf{u}}) = \sup \{ D(w) \bigm | w \in \mathcal{L}(\mathbf{u})\} \,.$$
\item Any infinite word with finite defect contains infinitely
many palindromes. \item Infinite words with zero defect correspond
exactly to rich words.
\end{enumerate}
Periodic words with finite defect have been studied in~\cite{BHNR}
and in~\cite{GlJuWiZa}. It holds that the defect of an infinite
periodic word with the minimal period $w$ is finite if and only if
$w=pq$, where both $p$ and $q$ are palindromes. Words with finite
defect have been studied in~\cite{BaPeSt} and~\cite{GlJuWiZa}.

\medskip

The number of  palindromes of a~fixed length occurring in an
infinite word is measured by the so called {\it palindromic
complexity} ${\mathcal{P}_{\bf u}}$, a~map which assigns to any
non-negative integer $n$ the number
$$ {\mathcal{P}_{\bf u}}(n) := \#\{ w \in {\mathcal{L}_n}(u)\mid w \ \
\hbox{is a palindrome}\}\,.$$

Denote by
$$
T_{\mathbf u}(n) = \mathcal{C}_{\mathbf
u}(n+1)-\mathcal{C}_{\mathbf u}(n)  +2 - \mathcal{P}_{\mathbf
u}(n+1) - \mathcal{P}_{\mathbf u}(n).
$$
The following proposition is proven in \cite{BaMaPe} for uniformly
recurrent words, however the uniform recurrence is not needed in
the proof, thus it holds for any infinite word with language
closed under reversal.

\begin{prop}[\cite{BaMaPe}]\label{Balazi}
 Let $\mathbf u$ be an infinite word
with language closed under reversal. Then
\begin{equation}\label{BalaziNerovnost}
T_{\mathbf u}(n)\geq 0\,,
\end{equation}
for all $n\in \mathbb{N}$.
\end{prop}

\medskip

It is shown in \cite{BuLuGlZa} that this bound can be used for
a~characterization of rich words as well. The following
proposition states this fact.

\begin{prop}[\cite{BuLuGlZa}]\label{rich_pal} An infinite word  ${\mathbf u}$
with language closed under reversal is rich if and only if the
equality in \eqref{BalaziNerovnost} holds for all $n\in
\mathbb{N}$.
\end{prop}

Let ${\mathbf u}$  be an infinite word with language closed under
reversal. Using the proof of Proposition~\ref{Balazi}, those $n \in
\mathbb N$ for which $T_{\mathbf u}(n)=0$ can be characterized in
the graph language.

An {\em $n$-simple path} $e$ is a~factor of ${\mathbf u}$ of
length at least $n + 1$ such that the only special (right or left)
factors of length $n$ occurring in $e$ are its prefix and suffix
of length $n$. If $w$ is the prefix of $e$ of length $n$ and $v$
is the suffix of $e$ of length $n$, we say that the $n$-simple
path $e$ starts in $w$ and ends in $v$. We will denote by
$G_n({\mathbf u})$ an undirected graph whose set of vertices is
formed by unordered pairs $(w,\overline{w})$ such that $w \in
\mathcal{L}_n({\mathbf u})$ is right or left special. We connect
two vertices $(w,\overline{w})$ and $(v,\overline{v})$ by an
unordered pair $(e,\overline{e})$ if $e$ or $\overline{e}$ is an
$n$-simple path starting in $w$ or $\overline{w}$ and ending in
$v$ or $\overline{v}$. Note that the graph $G_n({\mathbf u})$ may
have multiple edges and loops.

\begin{pozn}\label{graf_periodic}
Let us point out that if $\mathcal{L}_n({\mathbf u})$ contains no
special factor then $G_n({\mathbf u})$ is an empty graph. In this
case the word ${\mathbf u}$ is periodic, i.e., there exists a
primitive word $w$ such that ${\mathbf u} = w^\omega$ and $|w|\leq
n$. As proven in \cite{BHNR}, since the language of ${\mathbf u}$
is closed under reversal, the word $w$ is a product of two
palindromes. It is easy to see that $\mathcal{C}_{\mathbf
u}(n+1)=\mathcal{C}_{\mathbf u}(n)$ and $ 2=  \mathcal{P}_{\mathbf
u}(n+1) + \mathcal{P}_{\mathbf u}(n)$. Therefore $T_{\mathbf
u}(n)=0$.

\end{pozn}

\begin{lem} Let ${\mathbf u}$  be an infinite  word with language
closed under reversal, $n\in \mathbb{N}$.   Then $T_{\mathbf u}(n)
=
 0$ if and only if both of the following conditions are
met:  \begin{enumerate}\item The graph $G_n({\mathbf u})$  after
removing loops is a tree; \item Any $n$-simple path forming a loop
in the graph $G_n({\mathbf u})$  is a palindrome.
\end{enumerate}
\end{lem}
\pf It is a~direct consequence of the proof of Theorem 1.2 in~\cite{BaMaPe}
(recalled in this paper as Proposition~\ref{Balazi}). \pfk

\begin{coro}\label{podmnozina} Let ${\mathbf u}$  and ${\bf v}$ be  infinite
words with language closed under reversal  and $n \in \mathbb{N}$.
$$ \mathcal{L}_{n+1}({\bf v})\subset \mathcal{L}_{n+1}({\mathbf u}) \ \hbox{and } \  T_{\mathbf u}(n) =
 0 \qquad   \Longrightarrow \qquad T_{\bf v}(n) =
 0.$$
\end{coro}
\pf Our assumptions  imply that $G_n({\bf v})$ is a subgraph of
$G_n({\mathbf u})$ and  $G_n({\mathbf u})$ meets both conditions in
the previous lemma. These conditions are hereditary, i.e., any
connected  subgraph inherits these conditions as well. \pfk

\section{Brlek-Reutenauer conjecture}
Brlek and Reutenauer gave in~\cite{BrRe} a~conjecture relating the
defect and the factor and palindromic complexity of infinite words
with language closed under reversal.
\begin{conj}[Brlek-Reutenauer conjecture]\label{conjecture}
Let ${\mathbf u}$  be an infinite word with language closed under
reversal. Then
$$2D({\mathbf u}) = \sum_{n=0}^{+\infty}T_{\mathbf u}(n)\,.$$
\end{conj}
It is known from~\cite{BuLuGlZa} that the conjecture holds for
rich words.
\begin{thm}\label{rich_conj}
Conjecture~\ref{conjecture} is true if  ${\mathbf u}$  is rich.
\end{thm}
Brlek and Reutenauer provided in~\cite{BrRe} a~result for periodic
words.
\begin{thm}\label{periodic} Conjecture~\ref{conjecture} is true if  ${\mathbf u}$  is
periodic.
\end{thm}
In the sequel, we will prove the following theorem.
\begin{thm}\label{general}
Conjecture~\ref{conjecture} is true if $\mathbf u$ satisfies two assumptions:
\begin{enumerate}
\item Both $D(\mathbf u)$ and $ \sum_{n=0}^{+\infty}T_{\mathbf
u}(n)$ are finite.

\item For any $M\in \mathbb{N}$ there exists a factor $w \in
{\mathcal L}(\mathbf u)$ such that
\begin{itemize}
\item $w$ contains all factors of $\mathbf u$ of length $M$,
 \item $ww$ is a factor of $\mathbf u$.
\end{itemize}
\end{enumerate}
\end{thm}
In order to prove Theorem~\ref{general}, we need to put together
several claims. Let us first describe the main ideas of the proof.
The assumptions of Theorem~\ref{general} enable us to construct a
periodic word ${\mathbf v}$ with language closed under reversal
such that
\begin{itemize}
 \item   $D({ \mathbf u}) = D({ \mathbf v})$ and
 \item $T_{\mathbf v}(n)= T_{\mathbf u}(n)$ for all $n \in \mathbb{N}$.
\end{itemize}
Theorem~\ref{periodic} applied to the periodic word ${\mathbf v}$ then
concludes the proof.

Let us construct a~suitable periodic word. As $D(\mathbf u)$
is finite, there exists a factor $f\in {\mathcal{L}}({\mathbf{u}})$ such that $D(\mathbf u) = D(f)$.
Let us denote its length by $H=|f|$.
According to the inequality \eqref{BalaziNerovnost},  the finiteness of $
\sum_{n=0}^{+\infty}T_{\mathbf u}(n)$ implies that there exists an
integer $N\in \mathbb{N}$ such that $T_{\mathbf u}(n)=0$ for all
$n \geq N$. Let us put
\begin{equation}\label{maxM}
M=\max\{N,H\}.
\end{equation}
By Assumption (2), there
exists a  factor $w$ containing all elements of $ \mathcal{L}_M({
\mathbf u})$. Let us define
\begin{equation*}\label{def_v}
{\mathbf v} = w^\omega.
\end{equation*}

\begin{klejm} \label{claim1}
The word $w$ is a~concatenation of two palindromes, in particular,
the periodic word  $w^\omega$ has the language closed under reversal.
\end{klejm}

\pf Since the factor  $w$ contains the factor $f$ and the square
$ww$ belongs to ${\mathcal{L}}({\mathbf{u}})$, we have $D(f) \leq
D(w)\leq D(ww)\leq D({\mathbf{u}})$. As the factor $f$ was chosen
to satisfy $D(\mathbf u) = D(f)$, we may conclude that
\begin{equation}\label{vsechnorovno}
D(f) =D(w) = D(ww)= D({\mathbf{u}})\,.
 \end{equation}
The factor $ww$ is longer than the
factor $f$  and has the same defect as $f$.
Let us denote by $p$ the longest palindromic suffix  of $ww \in
\mathcal{L}({\mathbf u}) $. According to Corollary
\ref{prodlouzeni}, the palindrome $p$ occurs in $ww$ exactly once
and therefore $|p|
> |w|$. There exists a proper prefix $w'$ of  $w$ such that
$ww=w'p$. Let us denote by $w''$ the suffix of $w$ for which
$w=w'w''$. It means that $p=w''w'w''$. As $p$ is a palindrome, we
have $\overline{w''} = w''$ and $\overline{w'} = w'$. Hence the word $w$
is a~concatenation of two palindromes.

\pfk

\begin{klejm} \label{claim2}
$D({\mathbf v}) = D({\mathbf u})$\,.
\end{klejm}

\pf We will use  Theorem 6 from~\cite{BHNR}.  It implies that if
$w$ is product of two palindromes, then   $D(w^\omega) = D(ww)$.
This together with  \eqref{vsechnorovno} concludes the proof. \pfk

%
%

\begin{klejm} \label{claim3}
$T_{\mathbf v}(n)= T_{\mathbf u}(n)$ for all $n \in \mathbb{N}$.
\end{klejm}

\pf  Let us first consider $n \leq M-1$, where $M$ is the constant
given by~\eqref{maxM}. Since $w$ contains all elements of ${\mathcal L}_M(\mathbf u)$, we have $M \leq | w|$.
Since $ww \in \mathcal{L} ({\mathbf u})$, we also have $
\mathcal{L}_{M} ({\mathbf v}) =\mathcal{L}_{M} ({\mathbf u})$. It
implies
$$ \mathcal{C}_{\mathbf u}(n)=\mathcal{C}_{\mathbf v}(n) \quad \hbox{ and}\quad
\mathcal{P}_{\mathbf u}(n) = \mathcal{P}_{\mathbf v}(n)\quad
\hbox{ for all} \  n \leq M\,.$$ It gives the statement of the
claim for all $n \leq M-1$.

\medskip

Now we will consider $|w| > n \geq M$. According to the definition
of $N\leq M$, it holds $T_{\mathbf u}(n) = 0$. Since   $
\mathcal{L}_{n+1} ({\bf
 v})\subset
\mathcal{L}_{n+1}  ({\mathbf u}) $,  Corollary \ref{podmnozina}
gives $T_{\mathbf v}(n) = 0$ as well.

\medskip

Finally, we  consider $n \geq |w|\geq M$. Since $n$ is longer than or
equal to the period of ${\mathbf v}$ and since $w$ is a~product of
two palindromes, we have $\mathcal{C}_{\mathbf
v}(n+1)=\mathcal{C}_{\mathbf v}(n)$ and $\mathcal{P}_{\mathbf
v}(n+1)+\mathcal{P}_{\mathbf v}(n)=2$. It implies $T_{\mathbf
v}(n) = 0$. The value $T_{\mathbf u}(n)$ is zero as well,
according to the fact that  $N\leq M$. \pfk

\pf[Proof of Theorem~\ref{general}] It suffices to put together
Claims~\ref{claim1}, \ref{claim2}, and~\ref{claim3} and to realize
that Conjecture~\ref{conjecture} was already proven for periodic words. \pfk

\section{Brlek-Reutenauer conjecture holds for uniformly recurrent words}

In this section we will show that either both sides
in the Brlek-Reutenauer equality are infinite or   both
assumptions of Theorem~\ref{general} are satisfied for uniformly
recurrent words, which results in the main theorem of this paper.
\begin{thm}\label{hlavni} Conjecture~\ref{conjecture} is true if ${\mathbf u}$  is
uniformly recurrent.
\end{thm}

In order to prove Theorem~\ref{hlavni}, we will make use of
several equivalent characterizations of infinite words with finite
defect.
\begin{thm}\label{ekvivalence1}  Let ${\mathbf u}$  be a  uniformly recurrent  infinite
word with language closed under reversal. Then the following
statements are equivalent.
\begin{enumerate}
 \item The defect of ${ \mathbf u}$ is finite.
 \item  There exists an integer $K$ such that any complete return
 word of a~palindrome of length at least $K$ is palindrome as well.
\item There exists an integer $H$ such that the longest
palindromic suffix of any factor $w$ with length $|w| \geq H$
occurs in $w$ exactly once. \item There exists  an integer $N$
such that
$$T_{\mathbf u}(n) =
 0\quad \hbox{ for all} \ \ \   n\geq N\,.$$
\end{enumerate}
\end{thm}
\pf (1) and (2) are equivalent by Theorem 4.8 from~\cite{GlJuWiZa}. It
follows from the definition of $D(\mathbf u)$ that (1) and (3) are
equivalent. The equivalence of (1) and (4) was stated as Theorem
4.1 in~\cite{BaPeSt}. \pfk

\begin{coro}\label{nekonecno} Let ${\mathbf u}$  be a uniformly recurrent infinite
word with language closed under reversal. Then
$$\ \  D({\mathbf u}) \ \  \hbox{is finite}
\quad \Longleftrightarrow \quad  \sum_{n=0}^{+\infty}T_{\mathbf
u}(n)\ \ \hbox{is finite}\,.$$
\end{coro}

Thanks to Corollary~\ref{nekonecno}, we can focus on uniformly
recurrent words  ${ \mathbf u}$  with finite defect. An important
role in the proof of Theorem~\ref{hlavni} is the presence of
squares in ${\mathbf u}$.

\begin{lem}\label{ctverce} Let ${\mathbf u}$  be a uniformly recurrent infinite
word with finite defect and with language closed under reversal. Then
the set
$$ \{w \in \mathcal{A}^*\, |\,   ww \in \mathcal{L}({\mathbf u}) \}$$
is infinite.
\end{lem}
 \pf  We shall prove that for any $L\in \mathbb{N}$  there  exists
 a factor $w$  such that  $ww \in \mathcal{L}({\mathbf u}) $ and $|w| >L$.
WLOG take $L>K$, where $K$ is the constant from the statement (3) of Theorem~\ref{ekvivalence1}.
Then any complete return word of a~palindrome which is longer than $L$
is a palindrome as well. This implies that ${\mathbf u}$  has
infinitely many palindromes. Thus there exists an infinite
palindromic branch, i.e., a~both-sided infinite word $\ldots
v_3v_2v_1v_0v_1v_2v_3 \ldots $, where $v_i \in \mathcal{A}$ for
$i=1,2,3,\ldots$ and $v_0 \in \mathcal{A} \cup\{\epsilon\}$ such
that $v_k v_{k-1}\ldots v_0 \ldots v_{k-1}v_k\in
\mathcal{L}({\mathbf u})$ for any $k\in \mathbb{N}$. Consider a
palindrome $q=v_k v_{k-1}\ldots v_0 \ldots v_{k-1}v_k$ where $|q|>
3L$. Since ${\mathbf u}$ is uniformly recurrent, there exists an
index $i> k$ such that the factor $f=v_iv_{i-1}\ldots v_{k+2}v_{k+1}$ is
a~return word of $q$. The factor $fq$ is a~complete return word of
the palindrome $q$ and therefore $fq$ is a~palindrome.

At first suppose that the return word $f$ is longer then $|q|$. In
this case, $f=qp$ for some palindrome $p$. Hence the palindromic
branch has as its central factor the word $qpqpq$. We can put
$w=qp$.

Now suppose that the return word $f$ satisfies $|f|\leq |q|$. In
this case there exists an integer $j\geq 2$ and a factor $y$ such
that $fq = f^jy$  and $|y|< |f|$. If we  put $w=f^i$, with
$i=\lfloor \tfrac{j}{2}\rfloor$ then $ww \in \mathcal{L}({\mathbf
u})$ and $|w|> \tfrac{1}{3}|q|\geq L$. \pfk

\pf[Proof of Theorem~\ref{hlavni}] By Corollary~\ref{nekonecno},
the equality $2D({\mathbf u}) = \sum_{n=0}^{+\infty}T_{\mathbf
u}(n)$ holds as soon as one of the sides is infinite. Assume that
$D({\mathbf u})<+\infty$ and $\sum_{n=0}^{+\infty}T_{\mathbf
u}(n)<+\infty$. Let $M\in \mathbb{N}$ be an arbitrary integer. As
${\mathbf u}$ is uniformly recurrent, there exists an integer $R$
such that any factor longer than $R$ contains all factors  of
${\mathbf u}$  with length at most  $M$. According to Lemma
\ref{ctverce}, the set of squares occurring in ${\mathbf u}$ is
infinite, thus there exists a~factor $w$ longer then $R$ such that
$ww$ belongs to the language of ${\mathbf u}$. Its length
guarantees that $w$ contains all elements of
$\mathcal{L}_{M}(\mathbf u)$.

Consequently, Assumptions (1) and (2) of Theorem~\ref{general} are met
and the equality $2D({\mathbf
u}) = \sum_{n=0}^{+\infty}T_{\mathbf u}(n)$ follows.
\pfk

\section{Open problems}
In this section, we will summarize which statements concerning
defects are known for infinite words which are not necessarily
uniformly recurrent.

Let us transform Brlek-Reutenauer Conjecture into a~more
general question: ``For which infinite words $\mathbf u$ does the
equality
\begin{equation}\label{dulezite}
2D(\mathbf u)=\sum_{n=0}^{+\infty}T_{\mathbf u}(n)
\end{equation}
hold?''

In our summary of properties related to the above question,
let us first recall Proposition 4.6 from~\cite{GlJuWiZa} which
applies in full generality.
\begin{prop}\label{oddity} Let $\mathbf u$ be an infinite word.\\
$D(\mathbf u) \geq  \#\{\{v, \overline{v}\} \bigm |
v\not=\overline{v} \ \text{and} \ v \ \text{or} \ \overline{v} \
\text{is a complete return word in $\mathbf u$ of a palindrome
$w$}\}$.
\end{prop}
The set $\{v, \overline{v}\}$ is called an {\em oddity}.
\begin{obs}
If an infinite word $\mathbf u$ contains finitely many distinct
palindromes, then the equality~\eqref{dulezite} holds.
\end{obs}
\pf It follows from the definition that $D(\mathbf u)=+\infty$.
Since ${\mathcal P}_{\mathbf u}(n)=0$ for $n$ large enough and ${\mathcal C}_{\mathbf u}$ is non-decreasing, we have
$T_{\mathbf u}(n)\geq 2$ for such indices $n$. Consequently,
$\sum_{n=0}^{+\infty}T_{\mathbf u}(n)=+\infty$. \pfk
\begin{obs}
Let $\mathbf u$ be a~periodic word. Then the
equality~\eqref{dulezite} holds.
\end{obs}
\pf Theorem~\ref{periodic} states this fact  for infinite words with
language closed under reversal. In~\cite{BHNR} it is shown that
periodic words whose language is
not closed under reversal contain only finitely
many palindromes. Thus, the previous observation implies that the
equality is reached for such words, too. \pfk

From now on, let us limit our considerations to infinite words
containing infinitely many palindromes in their language.
\begin{obs} The equality~\eqref{dulezite} does not hold in general for infinite words which are not recurrent.
\end{obs}
\pf The word ${\mathbf u}=ab^\omega$ is rich, i.e., $D(\mathbf
u)=0$, however $\sum_{n=0}^{+\infty}T_{\mathbf u}(n)=-1$. \pfk

{\bf Problem 1}: It is an open problem whether the
equality~\eqref{dulezite} holds for recurrent words whose language
is not closed under reversal and contains infinitely many
palindromes. We have examples for which the equality holds and we
have so far no example refuting the equality~\eqref{dulezite}.
\begin{ex}
Let $\mathbf u$ be an infinite ternary word satisfying $\mathbf
u=\lim_{n \to +\infty}u_n$, where $u_0=a$ and
$u_{n+1}=u_nb^{n+1}c^{n+1}u_n$. The word $\mathbf u$ is recurrent,
however not closed under reversal (it does not contain the factor
$cb$). On one hand, $D(\mathbf u)=+\infty$ because $b^k$ has non-palindromic
complete return words for any $k \geq 1$, thus the number of
oddities is infinite. On the other hand,
since the only left extension of $a$ is $c$ and the only right extension of $a$ is $b$,
it is readily seen that the only palindromes of length $>1$
are of the form $b^n$ and $c^n$, thus ${\mathcal P}_{\mathbf u}(n)=2$ for all $n \geq 2$.
It is also easy to show that $c^n, \ b^n$, and $b^{n-1}c$ are distinct left special factors of length $n\geq 2$,
therefore ${\mathcal C}_{\mathbf u}(n+1)-{\mathcal C}_{\mathbf u}(n)\geq 3$ for all $n \geq 2$. This implies that $\sum_{n=0}^{+\infty}T_{\mathbf u}(n)=+\infty$.
\end{ex}
In the sequel, let us consider infinite words whose language is
closed under reversal and contains infinitely many palindromes.

Any rich word with language closed under reversal satisfies~\eqref{dulezite}
by Theorem~\ref{rich_conj}. For instance, the Rote word $\mathbf u$ - the fixed point of the morphism $\varphi$
defined by $\varphi(0)=001$ and $\varphi(1)=111$, i.e., $\mathbf u=\varphi(\mathbf u)$ - is rich because it satisfies $T_{\mathbf u}(n)=0$ for all $n \in \mathbb N$, which is not difficult to show. Therefore, the Rote word is an example of an infinite word which is not uniformly recurrent (it contains blocks of ones of any length) satisfying the equality~\eqref{dulezite}. We have, of course, no counterexample which would refute Brlek-Reutenauer Conjecture.

There exist several equivalent characterizations of words with
finite defect.
\begin{thm} \label{thm:cur_inf_eq}
Let $\mathbf u$ be an infinite word with language closed under
reversal and containing infinitely many palindromes. Then the
following statements are equivalent.
\begin{enumerate}
 \item The defect of ${ \mathbf u}$ is finite.
 \item $\mathbf u$ has only finitely many oddities.
\item There exists an integer $H$ such that the longest
palindromic suffix of any factor $w$ with length $|w| \geq H$
occurs in $w$ exactly once.
\end{enumerate}
\end{thm}
\pf (1) and (3) are equivalent by the definition of defect. (1)
implies (2) by Proposition~\ref{oddity}. The implication (2)
$\Rightarrow$ (1) was proved as Proposition 4.8 in~\cite{GlJuWiZa}
for uniformly recurrent words. However, we will show that the
proof works for words with language closed under reversal and
containing infinitely many palindromes, too.
\medskip

Assume that $D(\mathbf u)=+\infty$ and the number of oddities is
finite.

A finite number of oddities means that only finitely many
palindromes can have non-palindromic complete return words. Let
the longest such palindrome be of length $K$.

Since the number of palindromes is infinite, there exists an
infinite number of non-defective positions. Denote by $u^{(n)}$ the prefix of $\mathbf u$ of length $n$. Then $n$ is
a~non-defective position if $D(u^{(n-1)})=D(u^{(n)})$ (such positions correspond to the first occurrences
of palindromes).

There exists an integer $H$ such that the prefix of $\mathbf u$ of
length $H$ contains all palindromes of length $\leq K+2$. Hence,
if $n>H$ is a~non-defective position, then the longest palindromic
suffix of $u^{(n)}$ is of length greater than
$K+2$.

Since both the number of defective and non-defective positions is
infinite, we can find an index $k > H$ such that $k$ is
a~defective and $k+1$ a non-defective position. The longest palindromic
suffix $p$ of $u^{(k)}$ occurs at least twice in $u^{(k)}$, thus $u^{(k)}$ ends
in a non-palindromic complete return word of $p$. Since $k+1$ is a
non-defective position, it can be easily shown by contradiction
that the longest palindromic suffix of $u^{(k+1)}$ is of length
$\leq |p|+2 \leq K+2$.

This is a contradiction with the fact that non-defective positions
greater than $H$ have their longest palindromic suffix longer than
$K+2$. \pfk For words with language closed under reversal, some
implications remain valid. The first one is Proposition 4.3
and the second one is Proposition 4.5 from~\cite{BaPeSt}.
\begin{prop}
Let $\mathbf u$ be an infinite word with language closed under
reversal. Suppose that there exists an integer $N$ such that for
all $n\geq  N$ the equality $T_{\mathbf u}(n)=0$ holds. Then the
complete return words of any palindromic factor of length $n \geq
N$ are palindromes.
\end{prop}
\begin{prop}
Let $\mathbf u$ be an infinite word with language closed under
reversal. If there exists an integer $H$ such that for any factor
$f \in {\mathcal L}(\mathbf u)$ with $|f|\geq H$ the longest
palindromic suffix of $f$ is unioccurrent in $f$. Then $T_{\mathbf
u}(n)=0$ for any $n \geq H$.
\end{prop}
The last proposition together with Theorem \ref{thm:cur_inf_eq} results in the following corollary.
\begin{coro}
Let $\mathbf u$ be an infinite word with language closed under reversal with.
Then we have
$$
D({\mathbf u}) < + \infty \quad \Rightarrow \quad \sum_{n=0}^{+\infty}T_{\mathbf u}(n) < +\infty .
$$
\end{coro}
It is an open question whether the implications in the previous
propositions can be reversed.

{\bf Problem 2}: Let $\mathbf u$ be an infinite word with
language closed under reversal and containing infinitely many
palindromes. Assume that there exists an integer $K$ such that all
palindromes of length $\geq K$ have palindromic complete return
words. Does there exist an integer $N$ such that $T_{\mathbf
u}(n)=0$ for any $n \geq N$?

{\bf Problem 3}: Let $\mathbf u$ be an infinite word with
language closed under reversal and containing infinitely many
palindromes. Suppose that there exists an integer $N$ such that
for all $n\geq  N$ the equality $T_{\mathbf u}(n)=0$ holds. Does
there exist also an integer $H$ such that for any factor $f \in
{\mathcal L}(\mathbf u)$ with $|f|\geq H$ the longest palindromic
suffix of $f$ is unioccurrent in $f$?

We have seen that in the proof of the validity of Brlek-Reutenauer Conjecture for
uniformly recurrent words, an important role was played by the presence of big squares in such words.
This leads to the last open problem.

{\bf Problem 4}: Find other classes of infinite words containing for any $L$
a~factor $w$ such that $|w|>L$ and $ww$ belongs to the language.
\section{Acknowledgement}

We acknowledge financial support by the Czech Science Foundation
grant 201/09/0584, by the grants MSM6840770039 and LC06002 of the
Ministry of Education, Youth, and Sports of the Czech Republic,
and by the grant  of the Student Grant Agency of the Czech
Technical University in Prague.


\begin{thebibliography}{1}
\bibitem{BaMaPe} P. Bal\'a\v zi, Z. Mas\'akov\'a, E. Pelantov\'a, Factor versus palindromic complexity of uniformly recurrent infinite
words, Theoret. Comput. Sci. 380 (2007) 266–275.
\bibitem{BaPeSt}  L{\!'}. Balkov\'a, E. Pelantov\'a, \v S. Starosta, Infinite words with finite defect, Adv. in Appl. Math. (2011), doi: 10.1016/j.aam.2010.11.006, the original publication is available at http://www.sciencedirect.com/science/journal/01968858.
\bibitem{BHNR} S. Brlek, S. Hamel, M. Nivat, C. Reutenauer, On the palindromic complexity of infinite words, in: J. Berstel, J. Karhumäki, D. Perrin (Eds.), Combinatorics on Words with Applications, International Journal of Foundation of Computer Science, Vol. 15, No. 2, 2004, pp. 293-306.
\bibitem{BrRe} S. Brlek, C. Reutenauer, Complexity and palindromic defect of infinite words, Theoret. Comput. Sci. 412 Issues 4-5 (2011) 493-497.
\bibitem{BuLuGlZa} M. Bucci, A. De Luca, A. Glen, L. Q. Zamboni, A connection between palindromic and factor complexity
using return words, Adv. in Appl. Math 42 (2009) 60–74.
\bibitem{DrJuPi} X. Droubay, J. Justin, G. Pirillo, Episturmian words and some constructions of de Luca and
Rauzy, Theoret. Comput. Sci. 255 (2001) 539-553.
\bibitem{GlJuWiZa} A. Glen, J. Justin, S. Widmer, L. Q. Zamboni, Palindromic richness, Eur. J. Comb. 30 (2009) 510–531.
\end{thebibliography}
\end{document}